\documentclass[a4paper,12pt, reqno]{amsart}

\usepackage{amsmath, amsthm, amscd, amsfonts, amssymb, graphicx, color}
\usepackage[bookmarksnumbered, colorlinks, plainpages]{hyperref}

\usepackage{enumitem}
\usepackage{url}
\usepackage{multirow}
\usepackage{graphicx}
\usepackage{subfigure}
\usepackage[pass]{geometry}
\usepackage{cite}
\usepackage{cancel}
\usepackage{color}
\usepackage{pgf,tikz}
\usetikzlibrary{arrows}

\usepackage{lscape}
\tikzstyle{v} = [circle, draw, inner sep=2pt, minimum size=3pt, fill=black]
\usetikzlibrary{calc,shapes,decorations.pathreplacing,matrix}
\usetikzlibrary{patterns}

\tikzset{square matrix/.style={
    matrix of nodes,
    column sep=-\pgflinewidth, row sep=-\pgflinewidth,
    nodes={draw,
      minimum height=4.5pt,
      anchor=center,
      text width=4.5pt,
      align=center,
      inner sep=0pt
    },
  },
  square matrix/.default=1.2cm
}

\newtheorem{Theorem}{Theorem}[section]

\newtheorem{Definition}[Theorem]{Definition}

\newtheorem{Claim}[Theorem]{Claim}

\newtheorem{Corollary}[Theorem]{Corollary}
\newtheorem{Remark}[Theorem]{Remark}

\begin{document}
\pagestyle{empty}
\title{New classification of graphs in view of the domination number of central graphs}
\author[S. Fujita]{S. Fujita}
\address{Shinya Fujita, School of Data Science, Yokohama City University, 236-0027, Yokohama, Japan.}
\email{shinya.fujita.ph.d@gmail.com}
\author[F. Kazemnejad]{F. Kazemnejad}
\address{Farshad Kazemnejad, Faculty of Basic Sciences, Department of Mathematics, Ilam University, P. O.Box 69315-516, Ilam, Iran.}
\email{kazemnejad.farshad@gmail.com}
\author[B. Pahlavsay]{B. Pahlavsay}
\address{Behnaz Pahlavsay, Department of Mathematics, University of Mohaghegh Ardabili, 
P. O. Box 179, Ardabil, Iran.}
\email{pahlavsayb@gmail.com}

\date{\today}

\begin{abstract}
For a graph $G$, the central graph $C(G)$ is the graph constructed from $G$ by subdividing each edge of $G$ with one vertex and also by adding an edge to every pair of non-adjacent vertices in $G$. Also for a graph $G$, let $\gamma(G)$ and $\tau(G)$ be the domination number of $G$ and the minimum cardinarity of a vertex cover of $G$, respectively. 
In this paper, we give a new classification of graphs concerning the domination number of central graphs and minimum vertex covers of graphs. Namely, we show that any graph $G$ with at least three vertices can be classified into one of the two classes of graphs with $\gamma(C(G))=\tau(G)$ and $\gamma(C(G))=\tau(G)+1$, respectively, together with some special properties concerning a vertex cover of $G$. 
We also give some new results on the domination number of central graphs. 
\\[0.2em]

\\[0.2em]

\noindent 
\end{abstract}
\maketitle
\section{Introduction}


The notion of domination and its many generalizations have been intensively studied in graph theory and the literature on this subject is vast, see for example \cite{HHS5}, \cite{HHS6}, \cite{HeYe13},  \cite{dominmiddle},  \cite{3totdominrook} and \cite{dominLatin}. Throughout this paper, we use standard notation for graphs and we assume that each graph is finite, undirected and simple. For the
standard graph theory terminology not given here we refer to \cite{bondy2008graph}. Throughout this paper, $G$ is a non-empty, finite, undirected and simple graph with the vertex set $V(G)$ and the edge set $E(G)$.

Let $G$ be a graph with the vertex set $V(G)$ of \emph{order}
$n$ and the edge set $E(G)$ of \emph{size} $m$.
The \emph{open neighborhood} and the \emph{closed neighborhood} of a
vertex $v\in V(G)$ are $N_{G}(v)=\{u\in V(G)\ |\ uv\in E(G)\}$ and
$N_{G}[v]=N_{G}(v)\cup \{v\}$, respectively. 
We write $K_{n}$, $C_{n}$ and $P_{n}$ for a \emph{complete} graph, a \emph{cycle} graph and a \emph{path} graph of order $n$, respectively, while $G[S]$, $W_n$ and $K_{n_1,n_2,\cdots,n_p}$ denote the subgraph of $G$ induced on the vertex set $S$, a \emph{wheel} graph of order $n+1$, and a \emph{complete $p$-partite} graph, respectively. The \emph{complement} of a graph $G$, denoted by $\overline{G}$, is a graph with the vertex set $V(G)$ such that for every two vertices $v$ and $w$, $vw\in E(\overline{G})$\ if and only if $vw\not\in E(G)$. A \emph{vertex cover} of the graph $G$ is a set $D \subseteq V(G)$ such that every edge of $G$ is incident to at least one element of $D$. The \emph{vertex cover number} of $G$, denoted by $\tau(G)$, is the minimum cardinality of a vertex cover of $G$.
 An \emph{independent set} of $G$ is a subset of vertices of $G$, no two of which are adjacent. Also a \emph{maximum independent set} is an independent set of the largest cardinality in $G$. This cardinality is called the \emph{independence number} of $G$, and is denoted by $\alpha(G)$.

\vskip 0.2 true cm

Vernold et al., in \cite{Vv} introduced the following graph operation to construct a new graph called the \textit{central graph} from a given graph. 

\begin{Definition}
	\emph{\cite{Vv}} The central graph $C(G)$ of a graph $G$ of order $n$ and size $m$ is a graph of order $n+m$ and size $\binom{n}{2}+m $ which is obtained by subdividing each edge of $G$ exactly once and joining all the non-adjacent vertices of $G$ in $ C(G)$.
\end{Definition}

This notion attracts attention, and there are many works on central graphs. To name a few, Kazemnejad and Moradi \cite{Kaz19} investigated the total domination number of central graphs and Patil and Pandiya Raj \cite{PPR} explored the total graph of central graphs and their covering numbers. Jahfar and Chithra \cite{JC1, JC2} studied a further extension on central graphs by using some algebraic approach, and they established some applications for obtaining the number of spanning trees and the Kirchhoff index in those extended graphs. 

In this paper, we focus on the domination number of central graphs.

\vskip 0.2 true cm

\begin{Definition} \label{DS} A \emph{ dominating set}, briefly DS,  of a graph $G$ is a set $S\subseteq V(G)$
	such that  $N_G[v]\cap
	S\neq \emptyset$, for any vertex $v\in V(G)$. The \emph{ domination number of $G$ is the minimum cardinality of a DS of $G$ and is denoted by $\gamma(G)$}. Moreover, a  dominating set of $G$ of cardinality $\gamma(G)$ is called a $\gamma-$set of $G$.
\end{Definition}

The main purpose of this paper is to provide a new classification of general graphs $G$ in terms of the domination number of $C(G)$. 
We will show that any graph of order at least $3$ containing an edge belongs to one of three families of graphs; thus, it is a classification based on the properties to satisfy $\gamma(C(G))=\tau(G)$ or $\tau(G)+1$ (see Theorem~\ref{main} in Section 3).  

Let $G$ be a graph. For $S\subseteq V(G)$, let $N^*(S)=\{x\in V(G-S) ~|~S\subseteq N_G(x)\}$. By definition, note that $N^*(\emptyset)=V(G)$ and $N^*(\{x\})=N_G(x)$ for every $x\in V(G)$. 
In order to prove our main result, Theorem~\ref{main}, we introduce two integer valued functions $f(G)$ and $h(G)$ on $V(G)$ as follows.

\begin{center}

$h(G)=\min\{|S|+|E(G-S)|+i(N^*(S))~|~ \emptyset \neq S\subseteq V(G)\}$,
\end{center}
 where $i(N^*(S))=|\{ v\in V(G-S)|\ v\in N^*(S)$ and $v$ is an isolated vertex of $G-S\}|$. 

\begin{center}
$f(G)=\min\{|S|~|~  S\subseteq V(G) ~\text{and} ~S $ contains a dominating set of $ \overline{G}$ \\and a vertex cover of $G\}$.
\end{center}

As is shown in Section 2, in fact these two functions are ``nearly" equivalent (see Theorem~\ref{h=f}). 
To see this, we shall start with the following observation on $h(G)$ and $f(G)$:
 
\begin{Remark}\label{h(G)<f(G)}
For a graph $G$, we have $h(G)\leq f(G)\leq \tau(G)+1$.
\end{Remark}
\begin{proof}
Let $S$ be a subset of $V(G)$ such that $|S|=f(G)$ and $S$ contains a dominating set of $ \overline{G} $ and a vertex cover of $G$. Then we have $|S|=|S|+|E(G-S)|+i(N^*(S))\geq h(G)$.

To show the second inequality, let $C$ be a vertex cover of $G$ such that $|C|=\tau(G)$ and take a vertex $v\in V(G)\setminus C$. Clearly, $C\cup \{v\}$ forms a dominting set of $\overline{G}$. Thus, by definition, we have $f(G)\leq \tau(G)+1$.
\end{proof}
\begin{Remark}\label{exception}
For every pair $n, \ell$ with $n \geq \ell\geq 3$, $h(K_{\ell}\cup \overline{K_{n-\ell}})=\ell-1<f(K_{\ell}\cup \overline{K_{n-\ell}})=\ell$.
\end{Remark}
Note that, we here put $ K_{\ell}\cup \overline{K_{n-\ell}}=K_n $ for the case $n=\ell$ because we may regard 
$ \overline{K_0} $ as an empty graph with no vertices.

This paper organizes as follows. In Section 2, we discuss the equivalence of $h(G)$ and $f(G)$, and we will mainly show that the domination number of a central graph can be described as these two functions (see Theorems~\ref{h=f} and~\ref{gamma=h}). We also give some results on the domination number of a central graph. Section 3 provides our main result concerning a classification of graphs $G$ in terms of $\gamma(C(G))$ together with some further results. In Section 4, we provide other miscellaneous results on the domination number of central graphs.

\section{Results on $h(G)$ and $f(G)$ and its equivalence}

We firstly observe graphs $G$ with small $h(G)$. It is easy to characterize graphs $G$ such that $h(G)=1$.

\begin{Remark}
A graph $G$ of order $n\geq 1$ satisfies $h(G)=1$ if and only if $G=\overline{K_n}$.
\end{Remark}

A \textit{double star} is a tree of order at least $4$ with diameter exactly $3$. We have the following characterization for graphs $G$ such that $h(G)=2$.

\begin{Theorem}\label{h2}
A graph $G$ of order $n\geq 3$ such that $E(G)\neq \emptyset$ satisfies $h(G)=2$ if and only if either $G=K_3\cup \overline{K_{n-3}}$ or $G$ is a (possibly, disconnected) subgraph of a double star.
\end{Theorem}
\begin{proof}
It is easy to check that $(\Leftarrow)$ holds. So we shall prove $ (\Rightarrow) $ part. Since $h(G)=2$, note that $G$ contains at most two components of order at least $2$. Assume for the moment that $G$ contains such two components. Then obviously both components must be stars, and thus we can write $G=K_{1,\ell}\cup K_{1,m}\cup \overline{K_t}$, where $\ell\geq 1$, $m\geq 1$ and $t\geq 0$. Consequently $G$ is the subgraph of a double star, as desired. Thus we may assume that $G$ contains exactly one component $C$ such that $E(C)\neq \emptyset$. Assume for the moment that $C$ contains a cycle. Since $h(G)=2$, this forces $C\simeq K_3$, and it consequently follows that $G=K_3\cup \overline{K_{n-3}}$, as desired. Hence we may assume that $C$ is a tree such that $diam(C)\leq 4$. Suppose that $diam(C)=4$. Then $C$ contains an induced path $P=v_1v_2v_3v_4v_5$. Let $S$ be a subset of $V(G)$ such that  $|S|+|E(G-S)|+i(N^*(S))= h(G)=2$. Then we must have $S=\{v_2,v_4\}$ and $v_3\notin N^*(S)$. However, $v_3\notin N^*(S)$ implies that $E(G-S)\neq \emptyset$ and hence  $|S|+|E(G-S)|+i(N^*(S))>2$, a contradiction. So we have $diam(C)\leq 3$. This implies that $G$ is a desired graph because $G$ can be regarded as the subgraph of a double star.
\end{proof}

As observed from the proof of Theorem~\ref{h2}, in order to characterize graphs $G$ with $h(G)=c$ for a constant number $c$, we essentially need to investigate the structure of $G$ by looking at various subsets $S$ of $V(G)$ such that  $|S|+|E(G-S)|+i(N^*(S))=c$. This would require us a tedious argument when $c$ gets larger. 

We now turn our attention on the relationship between $h(G)$ and $f(G)$ for a given graph $G$.
In fact $h(G)=f(G)=\gamma(C(G))$ holds for almost all of graphs. Our results are following. 

\begin{Theorem}\label{h=f}
Let $G$ be a graph of order $n\geq 3$. If $G\neq K_{\ell}\cup \overline{K_{n-\ell}}$ 
for every $\ell$ with $n\geq \ell \geq 3$, then $h(G)=f(G)$.
\end{Theorem}
\begin{proof}
In view of Remark \ref{h(G)<f(G)}, it suffices to show that $ f(G)\leq h(G)$.
Let $S$ be a subset of $V(G)$ such that $|S|+|E(G-S)|+i(N^*(S))= h(G)$. Assume for the moment that $ i(N^*(S))>0 $ and let $v$ be an isolated vertex of $G-S$ such that $ N_G(v)\supseteq S$. Let $H$ be a vertex cover of $ G-S $ such that $|H|=\tau(G-S)$. Note that $|H|\leq |E(G-S)|$ (possibly $H=\emptyset$), and hence $|S\cup H\cup \{v\}|\leq h(G)$. Since $ S\cup H\cup \{v\} $ contains a dominating set of $\overline{G}$ and a vertex cover of $G$, we obtain $f(G)\leq h(G)$. Thus we may assume that $i(N^*(S))=0$. Again, let $H$ be a vertex cover of $G-S$ such that $|H|=\tau(G-S)$ (so $|H|\leq|E(G-S)|$). We may assume that $E(G-S)\neq \emptyset$ (otherwise, $|S|\geq f(G)$) and so $H\neq \emptyset$. If $S\cup H$ contains a dominating set of $\overline{G}$ and a vertex cover of $G$, we obtain $f(G)\leq h(G)$. Thus we may assume that there exists $v\in G\setminus (S\cup H)$ such that $N_G(v)\supseteq S\cup H$. Assume for the moment that $|S\cup H|=n-1$. In this case, since $G\neq K_n$, there exists $u\in V(G)$ such that $N_G(u)\nsupseteq V(G)\setminus\{u\}$. Since $V(G)\setminus\{u\}$ is a dominating set of $\overline{G}$ and a vertex cover of $G$, we obtain $f(G)\leq h(G)=n-1$. Thus we may assume that $|S\cup H|\leq n-2$ and there exists $v\in V(G)\setminus (S\cup H)$ such that  $N_G(v)\supseteq S\cup H$. Since $ i(N^*(S))=0 $ and $G\neq K_{\ell}\cup \overline{K_{n-\ell}}$ for any $n\geq \ell\geq 3$, there exist $z\in V(H)$ and $y\in V(G)\setminus (S\cup H\cup \{v\})$ such that $yz\in E(G)$ and $vy\notin E(G)$. Since $zv\in E(G)$, this implies that $|E(G-S)|\geq \tau(G-S)+1$. Note that $S\cup H\cup \{v\}$ is a dominating set of $\overline{G}$ and a vertex cover of $G$. Consequently, $f(G)\leq |S\cup H\cup \{v\}|\leq |S|+|E(G-S)|+i(N^*(S))= h(G)$. 
\end{proof}
\begin{Theorem}\label{gamma=h}
If $G$ is a graph of order $n\geq 3$ such that $G\neq K_{\ell}\cup \overline{K_{n-\ell}}$ 
for every $\ell$ with $n\geq \ell\geq 3$, then $\gamma(C(G))=h(G)$.
\end{Theorem}
\begin{proof}
We may assume that $\gamma(C(G))<f(G)$, and let $S$ be a $\gamma$-set of $C(G)$ (So $|S|<f(G)$).

\begin{Claim}\label{Clm1}
$S \nsubseteq V(G)$.
\end{Claim}
Suppose that $S\subseteq V(G)$. Then $S$ is a vertex cover of $G$ because $S$ is a $\gamma$-set of $C(G)$. Since $|S|<f(G)$, this implies that $S$ is not a dominating set of $\overline{G}$. Hence there exists $v\in V(G-S)$ such that $N_G(v)\supseteq S$. This means that $v$ is not dominated by $S$ in $C(G)$ because an internal vertex appears on every edge between $v$ and $S$ when we construct $C(G)$ from $G$. This is a contradiction, and thus the claim holds. \\ 

In view of Claim~\ref{Clm1}, $S$ can be partitioned into two part $S_1\neq \emptyset$ and $S_2\neq \emptyset$ such that $S_1\subsetneqq V(G)$ and $S_2\cap V(G)=\emptyset$. Note that $S_2$ forms a set of degree $2$ vertices in $C(G)$.

\begin{Claim}\label{Clm2}
$|S_2|\geq |E(G-S_1)|+i(N^*(S_1))$.
\end{Claim}
Let $x\in N^*(S_1)$ be an isolated vertex of $G-S_1$. Then there exists an edge $e=xy$ such that $y\in S_1$ in $G$. Since $S$ is a $\gamma$-set of $C(G)$, the vertex between $x$ and $y$ in $C(G)$ should be in $S_2$ (to dominate $x$). Also, since $V(G-S_1)\cap S_1=\emptyset$, any vertex $z$ appeared from an edge of $E(G-S_1)$ in $C(G)$ should be in $S_2$ (otherwise, $z$ is not dominated by $S$). This implies that the claim holds.\\

By Claim~\ref{Clm2}, we have $|S|\geq |S_1|+|E(G-S_1)|+i(N^*(S_1))\geq h(G)=f(G)$. Therefore, $\gamma(C(G))\geq h(G)$. Consequently, $\gamma(C(G))=h(G)$.
\end{proof}

Let $G$ be a graph. For a vertex cover $C$ of $G$ such that $|C|=\tau(G)$, we call $C$ a \textit{good cover} if $C$ is a dominating set of $\overline{G}$ (equivalently, $|C|=f(G)$). In other words, a good cover $C$ of $G$ is a minimum vertex cover of $G$ such that $V(G-C)\cap N^*(C)=\emptyset$ since otherwise $C$ cannot be a dominating set of $\overline{G}$. Note that the exceptional graphs $K_{\ell}\cup \overline{K_{n-\ell}}$ in Theorems~\ref{h=f} and~\ref{gamma=h}  contains no good cover. 

As an immediate corollary of Theorems~\ref{h=f} and~\ref{gamma=h}, we see that $\tau(G)$ can be a lower bound on $\gamma(C(G))$ for a graph $G$.  

\begin{Corollary}\label{lower bounds}
Let $G$ be a graph of order $n\geq 3$ such that $E(G)\neq\emptyset$. Then $\tau(G)\leq\gamma(C(G))$ and the equality holds only if $G$ contains a good cover. 
\end{Corollary}

We now point out that the independence number of a graph $G$ provides an upper bound on $\gamma(C(G))$. For a graph $G$ such that $E(G)\neq\emptyset$, let $S$ be an independent set of $V(G)$ (thus, $1\leq |S|\leq \alpha(G)<|V(G)| $). Let us observe that $G-S$ is a vertex cover of $G$ and $V(G-S)\cup \{x\}$ is a dominating set of $\overline{G}$ for any $x\in S$. 
This together with Remark~\ref{exception} and Theorems~\ref{h=f} and~\ref{gamma=h} implies the following corollary. 

\begin{Corollary}\label{independence number}
Any graph $G$ satisfies $\gamma(C(G))\leq \alpha(G)+1$ and the equality holds only if $S\cap N^*(G-S)\neq\emptyset$ for any maximum independent set $S$ of $G$.
\end{Corollary}

Utilizing Corollaries~\ref{lower bounds} and~\ref{independence number}, we obtain the following upper bounds on $\gamma(C(G))$ for bipartite graphs and disconnected graphs. 

\begin{Corollary}\label{bipartite graphs}
Let $G$ be a bipartite graph with partite sets $A$ and $B$ with $1\leq |A|\leq |B|$ such that $E(G)\neq\emptyset$. 
Then $\gamma(C(G))\leq |B|+1$ and the equality holds only if $B\cap N^*(A)\neq\emptyset$.
\end{Corollary}

\begin{Corollary}\label{theo:mindomincentraldis}
	Let $G$ be a graph of order $n\geq 3$ with no isolated vertex. If $G$ consists of a union of components $G_1, G_2, \cdots, G_\omega$ with $\omega\geq 2$, then
	\begin{center}
		$\tau(G_1)+\tau(G_2)+\cdots+\tau(G_\omega)\leq \gamma(C(G))\leq n-\omega$.
	\end{center}
\end{Corollary}

\begin{Remark}
There exist infinitely many graphs $G$ attaining the upper bound on $\gamma(C(G))$ in Corollaries~\ref{independence number} and~\ref{bipartite graphs}. 
To see this, for example, consider the balanced complete bipartite graph $K_{n,n}$ for $n\geq 2$ (see Theorem~\ref{theo:mindomincentralKn,m}). 
\end{Remark}

\section{Classification of graphs $G$ in terms of $\gamma(C(G))$}

In this section, we show that graphs can be classified into three types. More precisely, we have the following result. 
\begin{Theorem}\label{main}
Let $G$ be a graph with $n\geq 3$ vertices such that $E(G)\neq \emptyset$. Then one of the following holds:
\begin{enumerate}
\item  $G=K_{\ell}\cup \overline{K_{n-\ell}}$ for an integer $\ell$ with $n\geq \ell\geq 3$ (thus, $\gamma(C(G))=\tau(G)$).
\item $G$ contains a good cover and $\gamma(C(G))=\tau(G)$.
\item $G$ contains no good cover and $\gamma(C(G))=\tau(G)+1$.
\end{enumerate}
\end{Theorem}
\begin{proof}
Suppose that $G\neq K_{\ell}\cup \overline{K_{n-\ell}}$ for every $n\geq \ell\geq 3$.
If $G$ contains a good cover $C$, then it is easy to check that $C$ is a dominating set of $C(G)$ and $\gamma(C(G))=\tau(G)$.Thus (2) holds. 

Suppose that $G$ contains no good cover. 
By Remark~\ref{h(G)<f(G)}, this forces $f(G)=\tau(G)+1$. 
Hence, by Theorems~\ref{h=f} and~\ref{gamma=h}, we have $\gamma(C(G))=\tau(G)+1$. Thus (3) holds. 
\end{proof}
It is well known that computing $\tau(G)$ is NP-hard in general (see e.g.\cite{KV}). Therefore, by Theorem~\ref{main}, we have the following:
\begin{Corollary}
For a graph $G$, computing $\gamma(C(G))$ is NP-hard in general.
\end{Corollary}

We next consider so called ``vertex-preserving properties" in dominating sets of a graph. 
As can be seen in the construction of $C(G)$ from a graph $G$, there are a variety of graph operations in graph theory. More generally, by a graph operation $\mathcal{O}$, we often make a new graph $G'$ from a given graph $G$ (thus, we may put $G'=\mathcal{O}(G))$. When we consider a $\gamma$-set $S$ of $G'$, it often happens that $S$ cannot be constructed from the original vertices of $G$, meaning that $S\setminus V(G)\neq\emptyset$ (that is, $S$ contains a new vertex of $G'$ appeared from $G$ by the graph operation $\mathcal{O}$).

It would be interesting to consider a general question to ask what kind of graph operations $\mathcal{O}$ satisfy the property that there exists a $\gamma$-set of $\mathcal{O}(G)$ using the original vertices of $G$ only. The graph operation to construct the central graph from a graph satisfies this property to some degree. Namely we have the following result:

\begin{Theorem}
Let $G$ be a graph of order $n\geq 3$. Any $\gamma$-set $X$ of $C(G)$ satisfies $X\setminus V(G)\neq\emptyset$ if and only if $G\simeq K_{\ell}\cup \overline{K_{n-\ell}}$ for some $\ell$ with $n\geq \ell \geq 3$.
\end{Theorem}
\begin{proof}
It is easy to check that $(\Leftarrow)$ holds. So we prove $(\Rightarrow)$ part. Let $S$ be a vertex subset of $V(G)$ such that $|S|=f(G)$ and $S$ contains a dominating set of $\overline{G}$ and a vertex cover of $G$. Then it is obvious that $S$ can be a dominating set of $C(G)$. Then, in view of Theorems \ref{h=f} and \ref{gamma=h}, $S$ is a $\gamma$-set of $C(G)$; and moreover, since $S\setminus V(G)=\emptyset$, this forces $G\simeq K_{\ell}\cup \overline{K_{n-\ell}}$ for some $\ell$ with $n\geq \ell\geq 3$.
\end{proof}
\begin{Corollary}
For any graph $G$ of order $n\geq 3$ such that $G\neq K_{\ell}\cup \overline{K_{n-\ell}}$ for every $\ell$ with $n\geq \ell\geq 3$, there exists a $\gamma$-set $ D $ of $C(G)$ such that $D\subseteq V(G)$.
\end{Corollary}

\section{Domination numbers of central graphs}

Utilizing Theorem~\ref{main}, let us determine $\gamma(C(G))$ for various graphs $G$. We list the following results without proofs. One can easily check that these basic graphs $G$ in our theorems contain good covers and also it is easy to obtain $\tau(G)$, thereby proving that $\gamma(C(G))=\tau(G)$ by Theorem~\ref{main}. So we shall omit all the proofs.  

\begin{Theorem}\label{theo:mindomincentralspn}
\label{gamma_{t}(C(P_n))}
For any path $P_{n}$ of order $n\geq 3$, 
$$\gamma (C(P_n))=
\begin{cases}
 \lceil \frac{n}{2} \rceil & \text{ if $n=3,4,5,$}\\ 
\lfloor \frac{n}{2}\rfloor & \text{ otherwise.} 
\end{cases}
$$
\end{Theorem}

\begin{Theorem}\label{theo:mindomincentralC}
\label{gamma_{t}(C(C_n))}
For any cycle $C_{n}$ of order $n\geq 3$, 
$$\gamma (C(C_n))=
\begin{cases}
n-1 & \text{ if $n=3,4,$}\\ 
\lceil \frac{n}{2}\rceil & \text{ otherwise.} 
\end{cases}
$$
\end{Theorem}
\begin{Theorem}\label{theo:mindomincentralKn,m}
Let $K_{m,n}$ be the complete bipartite graph with $n\geq m \geq 2$. Then $\gamma(C(K_{m,n}))=m+1$. 
\end{Theorem}

\begin{Theorem}\label{theo:mindomincentralWn}
For any wheel $W_n$ of order $n+1\geq 4$, 
$$\gamma (C(W_n))=
\begin{cases}
n & \text{ if $n=4,$}\\ 
\lceil \frac{n}{2}\rceil+1 & \text{ otherwise.} 
\end{cases}
$$
\end{Theorem}

\begin{Definition} The \emph{friendship} graph $F_n$ of order $2n+1$ is obtained from $n$ copies of disjoint triangles $K_3$ by identifying one vertex of each triangle (as a common vertex).
\end{Definition}
\begin{Theorem}\label{theo:mindomincentralFn}
	Let $F_n$ be the friendship graph with $n\geq 2$. Then $\gamma(C(F_n))=n+1$. 
\end{Theorem}

We next consider a graph operation to construct a new graph $G'$ from more than one graphs and discuss the domination number of its central graphs. 

\begin{Definition}
The \emph{m-corona} $G\circ P_m$ of a graph $G$ is the graph of order $(m+1)|V(G)|$ obtained from $G$ by adding a path of order $m$ to each vertex of $G$. 
\end{Definition}

\begin{Theorem}\label{theo:mintotdomincorona}
For any connected graph $G$ of order $n\geq 3$, $$\gamma(C(G\circ P_1))=n.$$
\end{Theorem}
\begin{proof} 
By the construction, it is easy to see that $G$ forms a good cover of $G\circ P_1$. 
Thus, by Theorem~\ref{main}, we have $\gamma(C(G\circ P_1))=n.$
\end{proof}

We finally remark that $\gamma(\overline{C(G)})$ has a common small value for any graph $G$. 

\begin{Theorem}
	\label{gamma(overline{C(G)})=2}
	Let $G$ be a connected graph of order  $n \geq 4$. Then $\gamma(\overline{C(G)})=2$
\end{Theorem}
\begin{proof}
	Let $G $ be a connected graph of order $n \geq 4$ with the vertex set $V=\{v_{1},\ldots, v_{n}\}$. Then $V(C(G))=V(\overline {C(G)})=V\cup \mathcal{C}$ where $\mathcal{C}=\{ c_{i  j}:~v_{i}v_{j}\in E(G) \}$ and $ E(\overline {C(G)})=E(G)\cup\{c_{ij}v_{k}:~c_{ij}\in \mathcal{C},~v_k\in V, \mbox{ and }  k\neq i ,j\} \cup \{c_{ij} c_{i^{'}j^{'}}:~c_{ij}, c_{i^{'}j^{'}} \in \mathcal{C}\}$. Let $S$ be a $\gamma$-set of $\overline{C(G)}$. By Corollary~\ref{lower bounds}, one can easily check that $\gamma(\overline{C(G)})\geq 2$. If $G = K_{1, n-1} $, then $S=\{v_i , c_{jk}\}$ is a $\gamma$-set of $\overline {C(G)}$, we have $\gamma_{t}(\overline{C(G)})=2$. If $G \neq K_{1, n-1} $, then
	  there exist at least two edges $v_{i}v_{j}, v_{i^{'}}v_{j^{'}} \in E(G)$ such that $\{i,j\} \cap \{i^{'},j^{'}\}=\emptyset$. Since $S=\{c_{ij}, c_{i^{'}j^{'}}\}$ is a $\gamma$-set of $\overline {C(G)}$, we have $\gamma(\overline{C(G)})=2$.
\end{proof}

Combining Theorems~\ref{h2},~\ref{gamma=h} and~\ref{gamma(overline{C(G)})=2}, we obtain the following characterization for graphs $G$ to satisfy   $\gamma(C(G))=\gamma(\overline{C(G)})$.

\begin{Corollary}
Let $G$ be a connected graph of order  $n \geq 4$. 
Then the following statements are equivalent:
\begin{enumerate}
\item 
$\gamma(C(G))=\gamma(\overline{C(G)})$.
\item $\gamma(C(G))=2$.
\item
$G$ is either $G=K_3\cup \overline{K_{n-3}}$ or a subgraph of a double star containing at least one edge.
\end{enumerate}
\end{Corollary}

\textbf{Acknowledgment.}

Shinya Fujita's research was supported by JSPS KAKENHI (19K03603).


\end{document}